\documentclass{article}
\usepackage{amsmath,amsfonts}
\usepackage[english]{babel}
\def\xib{\mbox{\boldmath$\xi$}}

\def\Rbb{{\mathbb R}}
\def\E{{\bf E}}
\def\Eb{{\bf E}}
\def\C{{\bf C}}
\def\Hb{{\bf H}}
\def\P{{\bf P}}
\def\Pb{{\bf P}}
\def\n{\noindent}
\def\g{\gamma}
\def\a{\alpha}
\def\iy{\infty}
\def\l{\lambda}

\def\xb{{\bf x}}
\def\x{{\bf x}}
\def\Xb{{\bf X}}
\def\X{{\bf X}}
\def\Yb{{\bf Y}}
\def\s{\sigma}
\def\t{\theta}

\def\Cov{{\bf Cov}}
\def\Var{{\bf Var}}
\def\cof{\left(( 4/3)^m-(m/3)-1\right)}
\def\wx{\int_{I^m}\omega _{m}^{2}(\xb)d\xb}
\def\endsymbol{$\sqcup\mkern-12mu\sqcap$}
\def\done{\ \endsymbol\medskip}

\begin{document}

\title{On asymptotic efficiency of multivariate version of Spearman's rho}

\author{
Alexander Nazarov\footnote{Department of Mathematics and Mechanics, St.Petersburg State University,
28 Universitetskii Prospekt,
St. Petersburg, 198504, Russia} \
and \ Natalia Stepanova\footnote{School of Mathematics and Statistics,
Carleton University, 1125 Colonel By Drive, Ottawa, ON K1S 5B6, Canada}}

\date{}
\maketitle

\begin{abstract}
A multivariate version of Spearman's rho for testing independence is
considered. Its asymptotic efficiency is calculated under a general
distribution model specified by the dependence function. The
efficiency comparison study that involves other multivariate
Spearman-type test statistics is made. Conditions for Pitman
optimality of the test are established. Examples that illustrate the
quality of the multivariate Spearman's test are included.

\medskip


\end{abstract}

\section{Introduction}
Testing for independence among the components of
$m$-variate vector is an important statistical problem.
There is an extensive statistical literature on this topic.
Over the last two decades a variety of new multivariate measures of association have been suggested,
including those based on ranks, and their properties have been studied.

Let $\Xb_i=(X_{i1},\ldots,X_{im}),\;m\geq
2,\,i=1,\ldots,n,$ be independent random vectors with absolutely continuous
cdf $F$ and marginal cdfs $F_1,\ldots,F_m.$
 Denote by $R_{ij}$ the rank of $X_{ij}$ among
$X_{1j},\ldots, X_{nj},$ $i=1,\ldots,n,$ $j=1,\ldots, m.$
In case of bivariate random sample, when $m=2$, a commonly used statistic for
testing the hypothesis of independence, $H_0:F\equiv F_1 F_2$, is Spearman's correlation coefficient \cite{spear}
$$\rho_n=\frac{12}{n^2-1}\left\{n^{-1}\sum_{i=1}^n R_{i1}R_{i2}-\left(\frac{n+1}{2}\right)^2\right\}$$
that estimates the functional
\begin{gather}\label{sr}\rho(F)=12\int F dF_1 dF_2 -3.
\end{gather}

Among various multivariate extensions of Spearman's rho available in the statistical literature, the following three statistics seem to be quite popular
(see, for example, \cite{joe}, \cite{quessy}
\cite{s-s}, \cite{step}):
\begin{eqnarray}
S_{m,n}&=&\frac{1}{C_m}\left\{{n}^{-1}\sum_{i=1}^n \prod_{k=1}^m (n+1-R_{ij}) -
{\left(\frac{n+1}{2}\right)}^m\right\}\label{S}\\
W_{m,n}&=&\frac{1}{C_m}\left\{n^{-1}\sum_{i=1}^n \prod_{j=1}^m R_{ij} -
{\left(\frac{n+1}{2}\right)}^m\right\}\label{W},\\
V_{m,n}&=&\frac{12}{n^2 - 1} \left\{{m\choose
2}^{-1}\sum\limits_{1\leq j<j'\leq m}n^{-1} \sum_{i=1}^n
R_{ij}R_{ij'}-{\left(\frac{n+1}2\right)}^2\right\}.\label{V}
\end{eqnarray}
where $C_m={n}^{-1}\sum_{i=1}^n i^m - {\left({(n+1)}/{2}\right)}^m$
is a normalizing factor.

Statistic (\ref{V}) is simply the
average pair-wise Spearman's rho \cite[Ch. 6]{kend} that estimates \cite{joe}
\begin{equation}
\nu_m(F)=12\left\{{m\choose
2}^{-1}\int\sum_{j<j'}F_{j}F_{j'}dF\right\}-3,\nonumber
\end{equation}
Statistics (\ref{S}) and (\ref{W}) are natural generalization of Spearman's rho, as they are
sample counterparts of the functionals
\begin{eqnarray*}
s_m(F)&=&\frac{1}{d_{m}}\left\{\int F d F_1\ldots dF_m -c_{m}\right\},\\
w_m(F)&=&\frac{1}{d_{m}}\left\{\int F_1\ldots F_m dF-c_{m}\right\},
\end{eqnarray*}
 where $c_{m}=2^{-m},\;d_{m}=(m+1)^{-1}-2^{-m}$, respectively.
The correspondence between $S_{m,n}$, the main object under investigation in this paper, and $s_m(F)$ is easy to see. Indeed, let $F_n$ be
the multivariate empirical cdf that corresponds to $F$, and let
$F_{j,n}$ be the marginal empirical cdfs based on
$X_{1j},\ldots,X_{nj}$, $j=1,\ldots,m.$ Then
\begin{gather*}
R_{ij}=\sum_{k=1}^n \mathbb{I}\left(X_{kj}\leq X_{ij}\right)=nF_{j,n}(X_{ij})=(n+1)F^*_{j,n}(X_{ij}),
\end{gather*}
where $F^*_{j,n}=(n/(n+1))F_{j,n}$ are the modified empirical cdfs.
Therefore $S_{m,n}$ can be written in the form
\begin{gather*}
S_{m,n}=\frac{(n+1)^m}{C_m}\left(\int \prod_{j=1}^m(1-F^*_{j,n})dF_n -\frac{1}{2^m}\right),
\end{gather*}
where, taking into account that
$\sum_{i=1}^n i^m\sim {n^{m+1}}/{(m+1)}$, as $n\to \iy$,
we have
$${(n+1)^m}/{ C_m}\sim \left({1}/{(m+1)}-{1}/{2^m}\right)^{-1}={ 1}/{d_{m}}.$$
Thanks to the Glivenko--Cantelli theorem (see, for example, \cite[Sec. I.4, Th. 1]{borovkov})
the closeness of $S_{m,n}$ and $s_m(F)$ is now immediately seen by noting that
\begin{gather}\label{rel}
\int_{\Rbb^m}\prod_{j=1}^m
(1-F_j(x_j))dF(x_1,\ldots,x_m)=\int_{\Rbb^m} F(x_1,\ldots,x_m)
\prod_{j=1}^m dF_j(x_j).
\end{gather}
Equality (\ref{rel}) is easy to verify by integrating by parts on the left-hand side
and using the properties of a multivariate cdf.

All three measures of multivariate concordance, $\nu_m(F)$, $s_m(F)$, and $w_m(F)$, increase with respect to the {\it multivariate concordance ordering} introduced by Joe \cite[Sec. 2]{joe}, \cite[Ch. 2]{hjoe}. This ordering is based on the concept of {\it positive orthant dependence} \cite[Sec. 2.1]{hjoe}.
It results from a comparison of a multivariate random vector with a random vector of independent random variables having
the same univariate marginal distributions. More precisely, let $F$ and $G$ be two $m$-variate cdfs with corresponding survival functions
$\bar{F}$ and $\bar{G}$, i.e. $\bar{F}(x_2,\ldots,x_m)=\P_F(X_1>x_1,\ldots,X_m>x_m)$ and
$\bar{G}(x_2,\ldots,x_m)=\P_G(Y_1>x_1,\ldots,Y_m>x_m)$. Then $G$ is said to be {\it more concordant} than $F$  (written $F\prec_c G$)
if
$$F(\xb)\leq G(\xb)\quad\mbox{and}\quad \bar{F}(\xb)\leq \bar{G}(\xb),\quad \mbox{for all}\;\;\xb=(x_1\ldots,x_m)\in\Rbb^m.$$
That is, if $\Xb=(X_1,\ldots,X_m)\sim F$, $\Yb=(Y_1,\ldots,Y_m)\sim G$, and $F\prec_c G$, then the components of $\Yb$ are
more likely than those of $\Xb$ to take on small and large values simultaneously.
As shown in \cite{joe}, $F\prec_c G$ implies $\nu_m(F)\leq \nu_m(G)$
and $w_m(F)\leq w_m(G)$. The fact that $s_m(F)$ is also increasing with respect to $\prec_c$
follows immediately from
Lemma 3.3.1 of \cite{joe}, the Remark below this lemma, and equality
(\ref{rel}).

Unlike the classical problem of testing independence when $m=2$, there is still no
clear clear concept of {\it negative} multivariate concordance.
Some ``characterizations" of the negative multivariate concordance can be found, for example, in \cite{joe}.

In connection with testing independence among the components of a $m$-variate
random vector, statistics (\ref{S})--(\ref{V}) were studied by several authors. One of the earliest
comprehensive study related to multivariate rank statistics for testing independence can be found in \cite{ruym}.
The Pitman efficiency properties of the tests based on (\ref{S})--(\ref{V}) are investigated in \cite{joe}, \cite{quessy}, \cite{step}, among others.
The asymptotic normality of statistics (\ref{S})--(\ref{V}) is established in \cite{s-s}
under rather weak assumptions on the underlying distribution.
The asymptotic efficiency study
of the Spearman-type tests, including those based on $S_{m,n}$
and $V_{m,n}$, is conducted in \cite{quessy} under various distribution models.
The thorough study of Pitman efficiency properties of $W_{m,n}$ and $V_{m,n}$ is done in \cite{step}.

An interesting problem related to finding the Pitman efficiency of a
test is to discover the structure of the underlying distribution for
which the test is Pitman optimal. Many test statistics were
suggested by their authors empirically for solving particular
problems of testing hypotheses, and were supposed to work in one or
another particular situation. Problems of finding the most
favourable alternatives have been studied in \cite{behnen},
\cite{gregory}, \cite[Ch. 6]{nikbook}, \cite{nikpan}, \cite{step},
etc. In this paper,
assuming one-parameter model
\begin{gather}\label{model}
F_\theta({\bf x})=\prod\limits_{j=1}^{m}F_j(x_j)+\theta
\Omega_{m}(F_1(x_1),\ldots,F_m(x_m)),\quad{\bf
x}=(x_1,\ldots,x_m)\in \Rbb^m,
\end{gather}
where $\theta$ is a
parameter of association close to zero and $\Omega_m$ is the {\it
dependence function} defined on the unit $m$-cube and satisfying certain boundary and smoothness conditions,
we find the most favourable alternative to independence
for which the test based on $S_{m,n}$ is Pitman optimal.

In order to
determine the ``optimal" distribution function one has to solve a
variational problem of minimization of an appropriate functional
on a set of special type, depending on the structure of the test
statistic. Typically, optimality conditions for tests
are found by using the Lagrange multiplier rule applied to a
functional on a Banach space.
Under the validity of model (\ref{model}), the optimality problems
for the Spearman-type test statistics $W_{m,n}$ and $V_{m,n}$ have been solved in \cite{step}.
Compared to these two cases, the extreme problem
related to $S_{m,n}$ is much more complicated and is reduced to solving the system of partial
differential equations with non-standard boundary conditions.
In Section 4.2 we provide solution to a general $m$-dimensional extremal
problem that gives Pitman optimality conditions for the sequence  $\{S_{m,n}\}_{n\geq 1}$.

In Section 2 we introduce statistical model and describe its properties.
Some basic properties of the test statistic $S_{m,n}$, including its asymptotic normality in terms of the dependence function, are given in Section 3.
Asymptotic efficiency study is performed in Section 4. The key result of the paper, the Theorem of Section 4.3,
provides the most favourable alternative to independence for the test statistic at hand.

\section{Multivariate model}

\subsection{Definition of the model}
Suppose we observe an $m$-variate random sample ${\bf
X}_1,\ldots,{\bf X}_n$ of size $n$ from distribution $\Pb_\theta$ on
the measurable space $(\Rbb^m,{\cal B}^m)$ indexed by a parameter
$\theta\geq 0$. Then the full observation is a single observation
from the product $\Pb^{n}_{\theta}$ of $n$ copies of $\Pb_\theta$.
Let $F_\theta$ be the distribution function that corresponds to
$\Pb_\theta.$ When testing independence among the components of a
continuously distributed random vector, without loss of generality,
the marginal cdfs $F_j$, $j=1,\ldots, m,$ can be taken uniformly
distributed on the interval $[0,1]$. Then, the statistical model is
described in terms of distribution functions as the collection of
probability measures $\{\Pb^n_{\theta}:\theta\geq 0\}$ on the sample
space $(\Rbb^{m\times n},{\cal B}^{m\times n})$ such that
\begin{gather}\label{mod}
F_\theta({\bf x})=\prod\limits_{j=1}^{m}x_j+\theta
\Omega_{m}(\xb),\quad{\bf x}=(x_1,\ldots,x_m)\in [0,1]^m=I^m,\; m\geq 2,
\end{gather}
 is satisfied for sufficiently small value of $\theta$ subject to some restrictions on $\Omega_m$. To be precise,
let ${\cal{F}}_m=\{F_\theta\}$ be the class of absolutely
continuous cdfs of type (\ref{mod}) for which, cf. \cite[Sec. 2]{step},
\medskip

${\bf (C1)} \;\Omega_m({\bf x})\geq 0, \quad  \xb \in I^m,$

${\bf (C2)} \;\Omega_m({\bf x})|_{x_k=0}=0,\quad \Omega_m(1,\ldots,1, x_{k},1,\ldots,1)=0,\quad \; x_k\in[0,1],\quad 1\leq k\leq m$,

${\bf (C3)}\; \Omega_m({\bf x})|_{x_k=1}= \Omega_{m-1}(x_1,\ldots, x_{k-1}, x_{k+1},\ldots, x_m),\quad 1\leq k\leq m$,

${\bf (C4)}$ there exists a non-zero mixed derivative
\begin{equation}
\omega_m({\bf x})=\frac{\partial^{m}\Omega_m({\bf x})}{\partial
x_{i_{1}}\ldots
\partial x_{i_{m}}},\quad {\text {\rm for}}\,
\lambda_m{\text{\rm-almost all}}\; \xb\in I^m,\nonumber
\end{equation}
such that $\omega_m\in{\bf L}_2(I^m)$, where $\l_m$ is the Lebesgue measure on $(\Rbb^{m},{\cal B}^{m})$ and  $(i_1,\ldots,i_m)$ is
an arbitrary permutation of the set $\{1,\ldots,m\}.$

 Due to $(\bf C1)$, boundary conditions
${\bf (C2)}$, and the consistency property ${\bf (C3)},$ for
sufficiently small $\theta$, all the properties of a multivariate
cdf are satisfied. The regularity condition {\bf (C4)} implies local asymptotic normality of the sequence
of models $\{\Pb^n_{\t}:\t\geq 0\}$ at $\theta=0$ (see Section 2.2 for details).
In the sequel, $m$-variate sample $\X_1,\ldots,\X_n$ is assumed taken from distribution for which the cdf
$F_{\t}(\xb)$ belongs to ${\cal F}_m$, $m\geq 2.$ The symbols $\Eb_\theta$  and $\Var_\theta$
(with index $n$ omitted) are used below to denote the
 expectation and the variance with respect to $\Pb^n_\theta$.

 We are
interested in testing the hypothesis of independence
$$H_0:\;\theta=0$$ against the one-sided alternative
$$H_1:\theta>0. $$

In the case $m=2$, model (\ref{mod}) was first studied by Farlie
\cite{farlie} and appeared later in a number of publications (see
\cite[Sec. 1]{step} for references), sometimes with a specific
choice of dependence function. Considered under assumptions $({\bf
C1})$--$({\bf C4})$, model (\ref{mod}) is an extension of the
Farlie model to the multivariate case.

\subsection{Local asymptotic normality of the model} Recall that a
sequence of statistical models is {\it locally
asymptotically normal} (LAN) if it converges to a
Gaussian model whose properties are well known \cite[Sec. 2]{ih}, \cite[Sec. 7]{vaart}.

Let $f_\theta$ be the density of $\Pb_\theta$ with respect to
$\lambda_m,$ that is,
\begin{equation*}
f_\theta({\bf x})=1+\theta \omega_m({\bf
x}), \quad {\bf x}\in I^m,\quad m\geq 2,
\end{equation*}
and denote by $\dot{f}_\theta({\bf x})$ its partial derivative
with respect to $\theta$. The true statistical difficulty is to distinguish between
the null hypothesis and the alternative when
$\t$ is small, typically ``of size $O(n^{-1/2}).$''
Therefore we introduce a {\it local parameter}
$h=\sqrt{n}\,\theta$, and consider a local statistical experiment indexed by
$h$: $$(\Xb_1,\ldots,\Xb_n)\sim \{\Pb^n_{h/\sqrt{n}}:\,h\geq
0\}.$$ Our attention will be focused on the performance of the
test based on $S_{m,n}$ at alternatives $$H_{1n}:\,h>0$$
converging, as $n\to \infty,$ to the null hypothesis $$H_0:h=0.$$

Let $\Delta_{n,\theta}$ be a random vector such that
$\Delta_{n,\theta} \stackrel{{\rm d}}{\longrightarrow}{\cal
N}(0,I_{\theta})$, where for $\theta=\theta_n=h/\sqrt{n},$
$$I_{\theta}={\bf
E}_{\theta}\left(\frac{\partial}{\partial
\t}\log({d\Pb_\theta}/{d\lambda_m})\right)^2=\int_{I^m}\frac{\dot{f}^2_\theta({\bf
x})}{f_\theta(\xb)}\,d\xb,\quad \theta\geq0,$$ is the Fisher
information in the parametric family $\{f_\theta(\xb),\, \theta \geq
0\}.$ Thanks to Theorem 1.1 of \cite{ih}, under the regularity condition ${\bf (C4)}$, the sequence of
statistical experiments $\{\Pb^n_{h/\sqrt{n}}:h\geq 0\}$ is LAN at
the point $h=0$, that is, for any $h\geq 0$
\medskip
\begin{equation}\label{lan}
\log\frac{d\Pb^n_{h/\sqrt{n}}}{d\Pb^n_0}= h\Delta_{n,0}-\frac12
h^{2}I_0+o_{\Pb^n_0}(1), \quad n\to \iy.
\end{equation}
Under local asymptotic normality
\begin{gather*}
\log\frac{d\Pb^n_{h/\sqrt{n}}}{d\Pb^n_0}\stackrel{{\rm d}}{\longrightarrow}{\cal
N}\left(-\frac12 h^2I_0,h^2I_0\right),\quad n\to \iy,
\end{gather*}
and hence the sequences of distributions $\{\Pb^n_{h/\sqrt{n}}\}$ and
$\{\Pb^n_0\}$ are mutually contiguous (see \cite[Sec. 7.5]{vaart}). This
fact allows us to obtain, by means of Le Cam's third lemma \cite[Sec. 6.7]{vaart}, limit
distribution of $S_{m,n}$ under the sequence of alternatives
$H_{1n}$, once the limit distribution under $H_0$ is known. Another
useful consequence of the local asymptotic normality is the
existence of an upper bound on the asymptotic power function of the
test. This makes it possible to establish the conditions for
asymptotic optimality of the test statistic $S_{m,n}$ (see
\cite[Ch.15] {vaart}).

\section{Basic properties and asymptotic normality}
In this section we list some basic properties of the test statistic
$S_{m,n}$. First, note that $S_{m,n}$ is
symmetric in $m$ variables. It is normalized so that its value is 1
when $R_{i1}=R_{i2}=\ldots =R_{im}$, or equivalently,
$F_\t=\min(F_1,\ldots,F_m)$ (perfect positive dependence), and its
expected value under $H_0$ is zero. The lower bound of $s_m(F)$ is equal to
$$s_m\left(\max(F_1+\ldots +F_m-m+1,0)\right)=\frac{2^m(m+1)}{2^m-(m+1)}\left\{\frac{1}{(m+1)!}-\frac{1}{2^m}\right\},$$
which is $-1$ for $m=2$ and is greater than $-1$ for $m\geq 3$. The lower bound is an increasing function of $m$
tending to zero as $m$ gets larger.  Hence $S_{m,n}$, the sample version of $s_m(F)$,
also exceeds $-1$ and its lower bound tends to zero as $m$ increases.
For this reason, it is appropriate to use the statistic $S_{m,n}$ for testing the hypothesis independence $H_0:\theta=0$ against
the one-sided alternative $H_1:\theta>0$ only.
 The non-symmetry between the upper and lower bounds is due to the ``curse if
dimensionality" and is partly explained by the inequality
\cite[Lemma 3.8]{hjoe}
$$\max(F_1+\ldots +F_m-m+1,0)\leq F_{\t}\leq \min (F_1,\ldots,F_m),$$
where in contrast to the Fr\'{e}chet upper bound, $\min(F_1,\ldots,F_m)$, the Fr\'{e}chet lower bound, $\max(F_1+\ldots +F_m-m+1,0)$
is generally not a cdf, except for the case $m=2.$
Through the curse of dimensionality, the concepts of perfect positive and perfect negative dependence lose the symmetry of the
two-dimensional case.

There exist a variety of theorems on asymptotic normality of multivariate linear rank statistics.
A unifying approach to these various results is given, for example, in \cite{ruym}. In particular, Theorem 2 of \cite{ruym}
implies that $S_{m,n}$ is asymptotically normally distributed. For our purpose, however,
it is more convenient to establish asymptotic normality of $S_{m,n}$ through the
correspondence between $S_{m,n}$ and a closely related $U$-statistic.

The multivariate rank statistic $S_{m,n}$ is asymptotically equivalent to a $(m+1)$-dimensional
$U$-statistic, $U_{m,n}$, based on $\int F_\t dF_1\ldots dF_m.$
The kernel of the $U$-statistic comes from symmetrizing $\mathbb{I}(X_{m+1, j}<X_{jj}, j=1,\ldots,m),$ cf. \cite[eq. (3.4)]{step}:
\begin{gather*}
U_{m,n}={n\choose {m+1}}^{-1}\sum_{1\leq i_1<\ldots <i_{m+1}\leq n}g(\Xb_{i_1},\ldots,\Xb_{i_{m+1}}),
\end{gather*}
where
\begin{multline*}
g(\Xb_1,\ldots,\Xb_{m+1})=
\frac1{(m+1)!}\sum_{(i_1,\ldots,i_{m+1})}(\mathbb{I}(X_{i_{m+1},1}<X_{i_1,1},
\ldots,X_{i_{m},m}<X_{i_{m},m})-c_{m})/d_{m},
\end{multline*}
and the summation is extended over all permutations
$(i_1,\ldots,i_{m+1})$ of $\{1,\ldots, m+1\}.$

The following result establishes ``locally uniform" asymptotic normality of $U_{m,n}$.
It will be used for calculating the {\it slope} (or {\it efficacy}) of the test statistic $S_{m,n}$ whose
limit distribution coincides with that of $U_{m,n}$.

\n \textbf{Lemma 1.} \textsl{If $F\in{\cal F}_m=\{F_{h/\sqrt{n}}\}$, then}
for all $h\geq 0$,
$$\frac{\sqrt{n}(U_{m,n}-\mu_m(h))}{\s_m(h)} \stackrel{{\rm d}}{\longrightarrow}{\cal N}(0,1),\quad n\to \iy,$$
\textsl{where}
\begin{eqnarray*}
\mu_m(h)&=&\frac{2^m(m+1)}{2^m-(m+1)}\, h\int_{I^m}\Omega_m(\xb) d\xb,\\
 \s^2_m(h)&=&\s_m^2(0)=\frac{(m+1)^2}{(2^m -(m+1))^2}\left(\left(\frac43\right)^m-\frac{m}{3}-1\right).
\end{eqnarray*}
\bigskip

\n \textbf{Proof.} For $\t\geq 0$, put
$$\eta_m(\t)=\E_{\t}\Psi^2_{\t}(X_1)-\left(\E_{\t}U_{m,n}\right)^2,\quad\Psi_{\t}(\xb)=\E_{\t}g(\X_1,\ldots,\X_{m+1})|\X_1=\xb).$$
By the CLT for $U$-statistics (see, for example, \cite[Sec. 4.2]{korbor}) for all $\t\geq 0$
\begin{gather}\label{clt}
n^{1/2}((m+1)^2\eta_m(\t))^{-1/2}(U_{m,n}-\E_{\t}U_{m,n})\stackrel{{\rm d}}{\longrightarrow}{\cal
N}(0,1),\quad n\to \iy,
\end{gather}
provided $\eta_m(\t)>0$ and $\E_{\t}\Psi^2_{\t}(\X_1,\ldots,\X_{m+1})<\iy$.
Note that
\begin{eqnarray*}
\E_{\t} U_{m,n}&=&\E_{\t} g(\X_1,\ldots,\X_{m+1})\\&=&\E_{\t}\left(
\mathbb{I}(X_{{m+1},1}<X_{11},\ldots, X_{{m+1},m}<X_{m,m})-c_{m}\right)/d_{m}\\
&=&\left(\int_{\Rbb^m}F_{\t}(\xb)dF_1(x_1)\ldots dF_m(x_m) -c_{m}\right)\Big/d_{m}=s_m(F_{\t}).
\end{eqnarray*}
In particular, $\E_0 U_{m,n}=0.$
Next, under $H_0$
\begin{eqnarray*}
\Psi_0(\xb)&=&\frac{1}{(m+1)!}\sum_{(i_1,\ldots, i_{m+1})}\left\{\E_0\left(
\mathbb{I}(X_{{m+1},1}<X_{11},\ldots, X_{{m+1},m}<X_{m,m})|\X_1=\xb\right)-c_{m}\right\}/d_{m}\\
&=&\frac{1}{d_{m}}\left\{\frac{m!}{(m+1)!}\left(\frac{1}{2^{m-1}}\sum_{j=1}^m x_j+\prod_{j=1}^m (1-x_j)\right)-c_{m}\right\}\\
&=&\frac{1}{2^m -(m+1)}\left\{2\sum_{j=1}^m x_j+2^m\prod_{j=1}^m (1-x_j)-(m+1)\right\}.
\end{eqnarray*}
Under model (\ref{mod}), the calculation of
$\eta_m(0)=\E_0\Psi^2_0(\X_1)$ can be simplified by noting that in
case of independence, the vectors ${\bf 1}-\X_1$ and $\X_1$ are
equally distributed, each with i.i.d. uniform components. Therefore
\begin{eqnarray*}
\eta_m(0)&=&\E_0\Psi^2_0({\bf 1}-\X_1)=\frac{1}{{(2^m -(m+1))^{2}}}\,\E_0\left(2\sum_{j=1}^m (1-X_{1j})+2^m\prod_{j=1}^m X_{1j}-(m+1)\right)^2\\
&=&\frac{1}{(2^m -(m+1))^{2}}\left(\left(\frac43\right)^m-\frac{m}{3}-1\right).
\end{eqnarray*}
From this,
applying
(\ref{clt}) we get under $H_0$
\begin{gather*}
\sqrt{n}\sigma_m^{-1}(0)(U_{m,n}-\mu_m(0))\stackrel{{\rm d}}{\longrightarrow}{\cal N}(0,1),\quad n\to \iy,
\end{gather*}
where
\begin{gather}\label{ms}
\mu_m(0)=0,\quad\sigma_m^2(0)=\frac{(m+1)^2}{(2^m -(m+1))^2}\left(\left(\frac43\right)^m-\frac{m}{3}-1\right).
\end{gather}
Thus, for $\t=0$ the lemma is proved.

Now using the ``contiguity arguments" we will reduce the derivation of asymptotic normality under $\t_n=h/\sqrt{n}$
to derivation under $\t=0$.
First, applying the projection technique to the $U$-statistic $U_{m,n}$, we get
\begin{gather*}
\sqrt{n}(U_{m,n}-\mu_m(0))=\frac{1}{\sqrt{n}}\sum_{i=1}^n \psi_0(\X_i)+o_{\P_0}(1),
\end{gather*}
where
\begin{gather*}\psi_0(\xb)=(m+1)\Psi_0(\xb)=\frac{(m+1)}{2^m -(m+1)}\left\{2\sum_{j=1}^m x_j+2^m\prod_{j=1}^m (1-x_j)-(m+1)\right\}.
 \end{gather*}
 Then, under $H_{1n}: h>0$, Le Cam's third lemma implies (see \cite[Sec. 7.5]{vaart})
 \begin{gather*}
 \sqrt{n}(U_{m,n}-\mu_m(0))\stackrel{{\rm d}}{\longrightarrow}{\cal N}\left(h \E_0[\psi_0(\X_1)\dot{l}_0(\X_1)],\E_0\psi_0^2(\X_1)\right),
 \end{gather*}
where $\dot{l}_{\t}(\xb)=({\partial}/{{\partial\t}}) \log f_{\t}(\xb)={\omega_m(\xb)}/{(1+\theta\omega_m(\xb))},$ $\t\geq 0.$
In other words, the statistic $U_{m,n}$ is approximately normally distributed with
variance
$n^{-1}\E_0\psi_0^2(\X_1)=n^{-1}\s^2_m(0),$
where $\s_m^2(0)$ is defined in (\ref{ms}), and
mean value
\begin{eqnarray*}
\mu_m(h)&=&h \E_0[\psi_0(\X_1)\dot{l}_0(\X_1)]=h\int_{I^m}\psi_0(\xb)\omega_m(\xb) d\xb\\&=&\frac{(m+1) }{2^m-(m+1)}\,h
\int_{I^m}\left(2\sum_{j=1}^m x_j+2^m\prod_{j=1}^m (1-x_j)-(m+1)\right)\omega_m(\xb)d\xb.
\end{eqnarray*}
Notice that
\begin{gather*}
\int_{I^m}\omega_m(\xb)d\xb=0,\quad \int_{I^m}x_j\omega_m(\xb)
d\xb=0,\quad 1\leq j\leq m,\\
\int_{I^m}\prod_{j=1}^m (1-x_j) \omega_m(\xb) d\xb=\int_{I^m}\Omega_m(\xb)d\xb,
\end{gather*}
where the first two equalities are consequences of boundary conditions {\bf (C2)}, and the third one follows from (\ref{rel}).
Therefore
\begin{gather*}
\mu_m(h)=\frac{2^m(m+1)}{2^m-(m+1)}\, h\int_{I^m}\Omega_m(\xb) d\xb.
\end{gather*}
The proof is completed. \done
\medskip

Due to Lemma 1, the test based on $S_{m,n}$ rejects the null
hypothesis of independence at level approximately $\a$ if
$\sqrt{n}S_{m,n}/\s_m(0)>z_{\a}$, where $z_{\a}=\Phi^{-1}(1-\a)$ is
the quantile of order $(1-\a)$ of a standard normal distribution.

\section{Asymptotic efficiency} First, we calculate the Pitman efficiency of the test
statistic $S_{m,n}$. Denote by $\gamma_{m,n}(\t)$, $\t=h/\sqrt{n}\geq 0$, the power function of the
test of level approximately $\a$: $$\gamma_{m,n}(\t)={\bf
P}_\theta(\sqrt{n}S_{m,n}/\s_m(0)>z_{\a}).$$
If for a sequence of tests $\{T_n\}$ the corresponding sequence of power functions
satisfies $\gamma_{n}(h/\sqrt{n})\to 1-\Phi(z_\alpha
-h s)$, for every $h\geq 0$, then the sequence $\{T_n\}$ is said to have {\it slope} (or {\it efficacy}) $s$.
 A widely-recognized quantitative measure of
comparison of two statistical tests is the square of the quotient
of two slopes. This quantity is called the {\it asymptotic
relative efficiency} (ARE) of the tests. Further, if
the sequence of experiments $\{\Pb^n_\theta:\theta\geq 0\}$ is LAN
at $\theta=0$, then an upper bound on the slope exists
\cite[Th. 15.4]{vaart}. This yields the relative efficiency of
the test with slope $s$ and the best test and thus allows us to
determine the {\it absolute} quality of the former.

Lemma 1 implies that
the sequence $\{S_{m,n}\}_{n\geq 1}$ is locally uniformly asymptotically
normal. Then the general result on behavior of
the {\it local limiting} power function,
defined as $$\gamma_m(h)=\lim_{n\to \iy}\gamma_{m,n}(h/\sqrt{n}),\quad h\geq 0,$$
says that
$\gamma_m$ depends on the sequence ${\{S_{m,n}\}}_{n\geq 1}$ only through the
quantity $\mu_m^{\prime}(0)/\sigma_m(0)$, the slope of the sequence of tests (see
\cite[Th. 14.7]{vaart}).

\subsection{Relative and absolute measures of efficiency}
Next lemma gives an expression for the local limiting power function of the test at hand
in terms of the dependence function $\Omega_m$.

\bigskip

\n\textbf{Lemma 2.} \textsl{Assume model {\rm (\ref{mod})} and let $\g_m(h)=\lim\limits_{n\to \iy}\g_{m,n}\left({h}/{\sqrt{n}}\right).$ Then }
\begin{gather*}
\g_{m}\left({h}\right)=1-\Phi\left(z_{\a}-
\frac{2^m}{\left((4/3)^m-{m}/{3}-1\right)^{1/2}} \,h \int_{I^m}\Omega_m(\xb)d\xb\right)
\end{gather*}

\n \textbf{Proof.} In view of Lemma 1, the proof follows immediately from Theorem 14.7 of \cite{vaart}. \done

From Lemma 2, the measure of efficiency for the sequence $\{S_{m,n}\}_{n\geq 1}$ is equal to
\begin{gather}\label{slS}
\left(\frac{\mu_{m}^{\prime}(0)}{\s_{m}(0)}\right)^2=\frac{4^m}{(4/3)^m-m/3-1}\left(\int_{I^m}\Omega_m(\xb)d\xb\right)^2.
\end{gather}
For the multivariate Spearman-type statistics $W_{m,n}$ and $V_{m,n}$ these are (see \cite[Sec. 4.1]{step})
\begin{gather}\label{slW}
\left(\frac{\mu_{m,W}^{\prime}(0)}{\s_{m,W}(0)}\right)^2=\frac{4^m}{(4/3)^m-m/3-1}\left(\int_{I^m}\prod_{j=1}^m x_j \omega_m(\xb)\,d\xb\right)^2,
\end{gather}
and
\begin{gather}\label{slV}
\left(\frac{\mu_{m,V}^{\prime}(0)}{\s_{m,V}(0)}\right)^2=144{m\choose 2}^{-1}\left(\sum_{1\leq i<j\leq m}\int_{I_m}x_i x_j \omega_m(\xb)\,d\xb\right)^2,
\end{gather}
respectively. The asymptotic relative efficiency of $S_{m,n}$ relative to $W_{m,n}$ and $V_{m,n}$ is then
\begin{eqnarray*}
{\rm ARE}(S,W)&=&\left(\frac{{\mu_{m}^{\prime}(0)}/{\s_{m}(0)}}{{\mu_{m,W}^{\prime}(0)}/{\s_{m,W}(0)}}\right)^2,
\quad
{\rm ARE}(S,W)=\left(\frac{{\mu_{m}^{\prime}(0)}/{\s_{m}(0)}}{{\mu_{m,V}^{\prime}(0)}/{\s_{m,V}(0)}}\right)^2.
\end{eqnarray*}

\medskip
At this point, recall that the sequence of models $\{\P^n_{h/\sqrt{n}}:h\geq 0\}$ under consideration is
LAN at $h=0$. Therefore there exists an upper bound on the power function of the test (see \cite[Th. 15.4]{vaart}).
 More precisely, for all $h\geq 0$,
 \begin{gather*}
 \limsup_{n\to \iy}\g_{m,n}(h/\sqrt{n})\leq 1-\Phi(z_{\a}-h\sqrt{I_0}).
 \end{gather*}
That is, the square root of the Fisher information $I_0=\int_{I^m}\omega_m^2(\xb)\,d\xb$ is the largest possible slope:
\begin{gather*}
\left(\frac{\mu_m^{\prime}(0)}{\s_m(0)}\right)^2\leq \int_{I^m}\omega_m^2(\xb)\,d\xb,
\end{gather*}
or equivalently,
\begin{gather}\label{bound}
\frac{4^m}{(4/3)^m-m/3-1}\left(\int_{I^m}\Omega_m(\xb)d\xb\right)^2\leq \int_{I^m}\omega_m^2(\xb)\,d\xb.
\end{gather}
Therefore the Pitman absolute efficiency of the test based on $S_{m,n}$ is given by the formula
\begin{gather}\label{absS}
e_S(\Omega_m)=\frac{4^m}{\left((4/3)^m-m/3-1\right)}{\left(\int_{I^m}\Omega_m(\xb)d\xb\right)^2}\Big/{\int_{I^m}\omega_m^2(\xb)\,d\xb}.
\end{gather}
For a given function $\Omega_m$, the closer the value of $e_S(\Omega_m)$ to one, the better the test based on $S_{m,n}$.
Similarly, using (\ref{slW}) and (\ref{slV})
\begin{eqnarray}
e_{W}(\Omega_m)&=&\frac{4^m}{(4/3)^m-(m/3)-1}\left(\int_{I^m}
\prod\limits_{i}x_{i}\, \omega(\xb)d\xb\right)^2\Big/ \wx,\label{absW}\\
e_V(\Omega_m)&=& 144{m\choose
2}^{-1}\left(\sum_{i<j}\int_{I^m}x_{i}x_{j}\omega_m(\xb)d\xb\right)^2\Big/
\wx.\label{absV}
\end{eqnarray}

\subsection{Extremal problem} We are interested in finding the most favourable alternative,
determined by the dependence function $\Omega_m(\xb)$,
for which the sequence of test statistics $\{S_{m,n}\}_{n\geq 1}$ has the largest possible slope.
This problem is reduced to the problem of finding $\Omega_m(\xb)$ that delivers equality in inequality (\ref{bound}).
The latter is a particular
case of a general $m$-dimensional extremal problem treated below.

Let us introduce the space $\C_0^m$ of functions that are $m$-times continuously differentiable with respect to
each variable and obey certain boundary conditions:
$${\C}^m_0=\{\Omega\in\C^m(I^m):\Omega(\xb)|_{x_j=0}=0,\;j=1,\ldots,m\}.$$
Define a scalar product on $\C_0^m$ as follows:
\begin{gather}\label{sp}
(\Omega_1,\Omega_2)=\int_{I^m}\omega_1(\xb)\omega_2(\xb)d\xb,\quad \Omega_1,\Omega_2\in\C_0^m,
\end{gather}
where $\omega_i(\xb)=\dfrac{\partial^m \Omega_i(\xb)}{\partial
x_1\ldots \partial x_m},$ $i=1,2.$ Denote by $\Hb^m$ the closure of
the space $\C^m_0$ under the norm $\|\cdot\|$ induced by scalar
product (\ref{sp}). For any $m\geq 2$, $\Hb^m$ is a Hilbert space
whose properties are immediately derived from those for $m=2$
established in \cite{naznik}. In particular, the embedding of
$\Hb^m$ into $\C(I^m)$ is compact. Therefore, a function from
$\Hb^m$ equals zero on any ``left" side of the cube $I^m$ adjacent
to the origin.

Recalling condition {\bf (C2)} imposed on the dependence function $\Omega_m$, consider the problem of minimizing the functional
$\int_{I^m}\omega^2(\xb)d\xb$ on the subspace of $\Hb^m$ specified
by the boundary conditions on the ``right" sides of $I^m$ adjacent
to the point ${\bf 1}=(1,1,\ldots,1)$ provided $\int_{I^m}\Omega(\xb)d\mu(\xb)=1$,
with $\mu$ being a finite measure on $I^m$. In order to describe all
possible boundary conditions of this extremal problem we need some
notation.

Let $M=\{1,2,\ldots,m\}$ and let $2^M$ be the set of all subsets of $M$.
For any $U\subset M$, denote $\x_U$ the $|U|$-dimensional vector $\xb_U=(x_i:i\in U)$.
Then, any possible set of the boundary conditions has the form
\begin{gather*}
\Omega(\xb)|_{\xb_U={\bf 1}}=0,\quad U\in{\cal M},
\end{gather*}
where ${\cal
M}\subset 2^M$ is such that for any $U\subset V\subset 2^M$, $ U\in
{\cal M}$ implies $V\in{\cal M}.$ That is, if a set $U$ belongs to
${\cal M}$, then all its ``oversets" also belong to ${\cal M}.$ The
reason for this requirement is simple: if $\Omega\in {\Hb^m}$ takes
a zero value on the side $\{\xb_{U}={\bf 1}\}$, it also takes a zero
value on all the subedges of $I^m$ of less dimension.

\medskip
\n\textbf{Remark 1.} For any $U\in 2^M$ define an $m$-dimensional vector of Boolean variables $(y_j=\mathbb{I}(j\in U), j=1,\ldots,m)$.
Then $\mathbb{I}(U\in{\cal M})$ is a monotone Boolean function \cite{kl}. Denote by $N(m)$ the total number of such functions.
Obviously, the number of the above considered extremal problems is also equal to $N(m)$.
So far, no explicit formula for $N(m)$ as a function of $m$ has been found. For asymptotic behaviour of $N(m)$ as $m\to \iy$ see \cite{korsh}.

\medskip
 Return to the extremal problem of interest:
\begin{gather}\label{extpr}
\|\Omega\|_{\Hb^m}^2\to {\rm min},\quad\mbox{where}\quad \int_{I^m}\Omega(\xb) d\mu(\xb)=1,
\end{gather}
subject to the conditions
\begin{gather}\label{eq12}
 \Omega\in{\bf H}^m\quad\mbox{and}\quad
\Omega(\xb)|_{\xb_U=1}=0\quad\mbox{for all}\;\; U\in{\cal M}.
\end{gather}
For a set $U=(i_1,\ldots,i_l)\in{\cal M}$ and its complement (in
${M}$) $U^{c}=(j_1,\ldots,j_k),$ $l+k=m$, put
$$\xb_U\xb^2_{U^{c}}=x_{i_1}\ldots x_{i_l}x^2_{j_1}\ldots
x^2_{j_k},\quad \partial\xb_U \partial\xb^2_{U^{c}} =\partial
x_{i_1}\ldots \partial x_{i_l}\partial x^2_{j_1}\ldots \partial
x^2_{j_k},
$$
and define the functions
$$K_U(\xb,\xib)=K_{i_1}(\xb,\xib)\ldots K_{i_l}(\xb,\xib),\quad k_{U^c}(\xb,\xib)=k_{j_1}(\xb,\xib)\ldots k_{j_k}(\xb,\xib),\quad\xb, \xib\in I^m,$$
where $$K_j(\x,\xib)=\min(x_j, \xi_j ),\quad k_j(\xb,\xib)=x_j \xi_j,\quad j=1,\ldots,m.$$
According to the Lagrange principle applied to a functional on a Banach space (see \cite[Sec. 2.2.3]{atf}),
the necessary condition of a minimum in (\ref{extpr})--(\ref{eq12}) is reduced to the Euler--Lagrange equation
\begin{gather}\label{bp1}
(-1)^m\lambda\frac{\partial^{2m}\Omega(\xb)}{\partial x_1^2\ldots \partial
x_m^2}=\mu(\xb),
\end{gather}
and the natural boundary conditions
\begin{gather}\label{bp2}
\frac{\partial^{l+2k}\Omega(\xb)}{\partial{\xb_V} \partial\xb^2_{V^c}(\xb)}\Big|_{\xb_V =1}=0,\quad\mbox{for any}\;\;
V\notin{\cal M},\; V\neq \emptyset,
\end{gather}
where the Lagrange multiplier $\l$ is found from the integral restriction in (\ref{extpr}).
The following result holds true.

\bigskip
\n\textbf{Lemma 3.} \textsl{Solution to extremal problem {\rm (\ref{extpr})}--{\rm(\ref{eq12})} is given by the formula
\begin{gather*}
\Omega(\xb)=\l^{-1}\int_{I^m}{\cal G}_{\cal M}(\xb,\xib)d\mu(\xib),\quad \xb\in I^m,
\end{gather*}
where ${\cal G}_{\cal M}$ is the Green function of boundary-value problem {\rm(\ref{eq12})}--{\rm (\ref{bp2})}
equal to
\begin{gather}\label{Gf}
{\cal G}_{\cal M}(\xb,\xib)=K_M(\xb,\xib)-\sum_{U\in{\cal M}}a_U K_{U^c}(\xb,\xib) k_{U}(\xb,\xib),
\end{gather}
with the coefficients $a_U$ defined recurrently by
\begin{gather}\label{rec}
\sum_{V\subset U\atop V\in{\cal M}}a_V=1,\quad\mbox{for all}\; U\in{\cal M},
\end{gather}
 and the constant $\l$ is given by
\begin{gather}\label{lam}
\l=\iint_{I^m\times I^m}{\cal G}_{\cal M}(\xb,\xib)d\mu(\xb)d\mu(\xib).
\end{gather}
}
\bigskip

\n\textbf{Proof.} First, note that
\begin{gather*}
-\frac{\partial^2 K_j(\xb,\xib)}{\partial{x_j^2}}=\delta(x_j-\xi_j),\quad \frac{\partial^2 k_j(\xb,\xib)}{\partial{x_j^2}}=0,
\end{gather*}
where $\delta$ is the Dirac function. Therefore the function ${\cal G}_{\cal M}$ in (\ref{Gf}) satisfies
$$(-1)^m \frac{\partial^{2m} {\cal G}_{\cal M}(\xb,\xib)}{\partial{x_1^2}\ldots \partial x_m^2}=\delta(\xb-\xib),$$
with an arbitrary choice of the constants $a_U$. The function ${\cal G}_{\cal M}$ also satisfies natural boundary conditions
(\ref{bp2}).

Taking into account (\ref{eq12}), we arrive at recurrent system (\ref{rec}).
Thus, solution to boundary-value problem (\ref{eq12})--(\ref{bp2}) is given by (\ref{Gf}).

It remains to note that the Lagrange multiplier $\lambda$ is found from the integral restriction in (\ref{extpr}) and
has the form (\ref{lam}). The lemma is proved. \done

\n\textbf{Remark 2.} Consider the following three sets of boundary conditions:
(i) there are no restrictions on $\Omega\in{\Hb}^m$ except for those that specify the space ${\Hb}^m$.
(ii) $\Omega\in{\Hb}^m$ equals zero on any $(m-1)$-dimensional side of $I^m$, and
(iii) $\Omega\in{\Hb}^m$ equals zero at the point ${\bf 1}=(1,\ldots, 1)$.
Then ${\cal M}=\emptyset$, ${\cal M}=2^M$, and ${\cal M}=\{M\}$,
respectively, and by Lemma 3 the corresponding Green functions are
\begin{eqnarray*}
{\cal G}_{\emptyset}(\xb,\xib)&=&\prod_{j=1}^m K_j(\xb,\xib),\\
{\cal G}_{2^M}(\xb,\xib)&=&\prod_{j=1}^m (K_j(\xb,\xib)-k_j(\xb,\xib)),\\
{\cal G}_{M}(\xb,\xib)&=&\prod_{j=1}^m K_j(\xb,\xib)-\prod_{j=1}^m k_j(\xb,\xib).
\end{eqnarray*}
These are covariance functions of the classical Gaussian random fields. They
correspond to a Brownian sheet, a Brownian pillow,
and a ``pinned" Brownian sheet, respectively, that
emerge as limiting processes in
nonparametric testing of multivariate independence.
For example, in the case $m=2$, the functions ${\cal G}_{2^M}$ and ${\cal G}_{M}$ appeared
in connection with finding the approximate Bahadur efficiency of independence tests
based on the comparison of the multivariate empirical cdf $F_n$ with
the product of margins $\prod_{j=1}^m F_j$ and with the product of empirical margins
$\prod_{j=1}^m F_{j,n}$ (see \cite[Ch. 5]{nikbook} for details).

\subsection{Most favourable alternative to independence} In order to determine
the ``optimal" distribution function for the sequence $\{S_{m,n}\}_{n\geq 1}$ we have to solve a variational problem
of minimization of the functional $\int_{I^m}\omega_m^2(\xb)d\xb$ on a set of functions of special type (see inequality (\ref{bound})).
Optimality conditions for the test statistic $S_{m,n}$ are given by the following theorem.

\bigskip
\n\textbf{Theorem.} \textsl{Let $F_{\t}\in{\cal F}_m$. Then the sequence of test statistics $\{S_{m,n}\}_{n\geq 1}$
is Pitman optimal if and only if}
\begin{gather}\label{op}
\Omega_m(\xb)=C\prod_{j=1}^m x_j \left(\prod_{j=1}^m (2-x_j)+\sum_{j=1}^m x_j -(m+1)\right),\\
\xb=(x_1,\ldots,x_m)\in I^m,\quad C>0.\nonumber
\end{gather}

\bigskip
\n\textbf{Proof.}
The test based on $S_{m,n}$ is the ``best"
for those dependence functions $\Omega_m$ that deliver equality in inequality (\ref{bound}).
Thus, we minimize the functional $\int_{I^m}\omega^2_m(\xb)d\xb$
on the space $\Hb^m$ subject to
\begin{gather*}
\int_{I^m}\Omega_m(\xb)d\xb =1,\quad \Omega_m(\xb)|_{x_{i_1}=\ldots= x_{i_{m-1}}=1}=0,\quad 1\leq i_1<\ldots <{i_{m-1}}\leq m,
\end{gather*}
where the second constraint on $\Omega_m$ is a consequence of
condition {\bf (C2)}. Therefore, with the notation of Section 4.2
$${\cal M}=\{M,M\setminus\{1\},M\setminus\{2\},\ldots,M \setminus\{m\}\},$$ and the problem (\ref{eq12})--(\ref{bp2})
takes the form
\begin{gather*}
(-1)^m\l\frac{\partial^{2m}\Omega_m(\xb)}{\partial x_1^2\ldots\partial x_m^2}=1,\\
\int_{I^m}\Omega_m(\xb)d\xb=1,\quad \Omega_m(\xb)|_{x_j=0} =0,\quad j=1,\ldots,m,\\
\Omega_m(\xb)|_{x_{i_1}=\ldots= x_{i_{m-1}}=1}=0,\quad 1\leq i_1<\ldots <{i_{m-1}}\leq m,\\
\dfrac{\partial^{2m-1}\Omega_m(\xb)}{\partial x_{i_1}\partial x_{i_2}^2\ldots \partial x_{i_m}^2}\,\Big|_{x_{i_1}=1}=0,\quad 1\leq i_1\leq m,\\
\dfrac{\partial^{2m-2}\Omega_m(\xb)}{\partial x_{i_1}\partial x_{i_2}\partial x_{i_3}^2\ldots \partial x_{i_m}^2}\,\Big|_{x_{i_1}=x_{i_2}=1}=0,\
\quad 1\leq i_1<i_2\leq m,
\\ \vdots\\
\dfrac{\partial^{m+2}\Omega_m(\xb)}{\partial x_{i_1}\ldots \partial x_{i_{m-2}}\partial x_{i_{m-1}}^2 \partial x_{i_m}^2}\,\Big|_{x_{i_1}=\ldots =x_{i_{m-2}}=1}=0,\quad 1\leq i_1<\ldots <i_{m-2}\leq m.
\end{gather*}

\n According to Lemma 3 the minimum of $\int_{I^m}\omega_m^2(\xb) d\xb$
is attained for the function
\begin{gather}\label{opt}
\Omega_m(\xb)=\l^{-1}\int_{I^m}{\cal G}(\xb,\xib) d\xib,
\end{gather}
where
$${\cal G}(\xb,\xib)=\prod_{j=1}^m K_j(\xb,\xib)-\sum_{j=1}^m \left(K_j(\xb,\xib)\prod_{i\neq j} k_i(\xb,\xib)\right) +(m-1)\prod_{j=1}^m k_j(\xb,\xib),$$
with $K_j$ and $k_j$ as before. By homogeneity of inequality (\ref{bound}) the extremal function is defined up to a positive constant.
Integrating in (\ref{opt}) yields (\ref{op}). \done

\n\textbf{Remark 3.} For the Spearman-type test statistics $W_{m,n}$
and $V_{m,n}$ the most favourable alternatives are specified by the
dependence functions (see \cite[Sec. 5]{step})
\begin{eqnarray*}
\Omega_{m,W}(\xb)&=&C\prod\limits_{j=1}^{m}x_j
\left(\prod\limits_{j}x_j-\sum\limits_{j}x_j+(m-1)\right),\quad \xb\in I^m,\quad C>0,\\
\Omega_{m,V}(\xb)&=&C\prod\limits_{i=j}^{m}x_j \sum\limits_{i<j}(1-x_i)(1-x_j),\quad \xb\in I^m,\quad C>0,
\end{eqnarray*}
respectively. The function $\Omega_{m,V}$, that corresponds to
the pair-wise average Spearman's statistic $V_{m,n}$, determines an $m$-variate extension of the Farlie--Gumbel--Morgenstern
distribution introduced in \cite[Sec. 5.1]{hjoe}.

\subsection{Examples} Now we examine, for several choices of cdf $F_\t$,
the Pitman efficiency of $S_{m,n}$ compared to the other two multivariate
Spearman-type test statistics, $W_{m,n}$ and $V_{m,n}$.

\medskip
\n {\bf Example 1.} Let $\Xb_1,\ldots,\Xb_n$ be independent copies of
the equicorrelated random Gaussian vector $\Xb=(X_1,\ldots,X_m)$
with $\Eb X_i=0,$ $\Var X_i=1$, and $\Cov(X_i,X_j)=\theta,$ $1\leq i\neq j\leq m,$
so that the experiment $\{\P_{\t}^n:\t\geq 0\}$ is normal.
The first two terms of
Taylor's expansion of the cdf of $\Xb$ around
$\theta$ are of the form (\ref{mod}) with the dependence function \cite{genest}, \cite{step}
\begin{equation}
\Omega_{m}(\xb)=\sum_{1\leq i<j\leq
m}\varphi(\Phi^{-1}(x_{i}))\varphi(\Phi^{-1}(x_{j}))\prod_{k\neq
i,j} x_k,\quad \xb\in I^m.\nonumber
\end{equation}
The mixed derivative of $\Omega_{m}(\xb)$ is
$$\omega_{m}(\xb)=\sum\limits_{1\leq i<j\leq
m}\Phi^{-1}(x_{i})\Phi^{-1}(x_{j}),$$ and
\begin{eqnarray*}
\int_{I^m}x_{i}x_{j}\,\omega_{m}(\xb)d\xb
&=&{1}/({4\pi}),\quad \int_{I^m}\omega^2_{m}(\xb) d\xb ={m(m-1)}/{2},\\
\quad \int_{I^m}\prod\limits_{i}x_i \,\omega_{m}(\xb)d\xb&=&\int_{I^m}\Omega_m(\xb) d\xb={m(m-1)}/({2^{m+1}\pi}).
\end{eqnarray*}
Now applying (\ref{absS})--(\ref{absV}) we obtain
\begin{equation*}
e_S(\Omega_m)=e_{W}(\Omega_{m})=\frac{m(m-1)}{2\pi^2\cof},\quad e_{V}(\Omega_{m})=\frac{9}{\pi^2}\approx 0.9119;
\end{equation*}

The asymptotic efficiency of $S_{m,n}$ and $W_{m,n}$ decreases in $m$ and
equals 0.8207, 0.7349, 0.6548 for $m=3,\,4,\,5,$ respectively,
whereas the asymptotic efficiency of $V_{m,n}$ is a
constant close to 1 and independent of $m$. Thus, in the normal
case, the average test based on $V_{m,n}$ is asymptotically more efficient than the
multivariate Spearman's tests based on $S_{m,n}$ and $W_{m,n}$.

\medskip
\n \textbf{Example 2.} Consider the multivariate extension of the Farlie--Gumbel--Morgenstern distribution for which
the dependence function is
$$\Omega_{m}(\xb)=\prod\limits_{i=j}^{m}x_j \sum\limits_{i<j}(1-x_i)(1-x_j),\quad \xb\in I^m.$$
In this case, the average pair-wise Spearman's test based on $V_{m,n}$ is Pitman optimal, i.e.,
$e_V(\Omega_m)=1$ (see \cite[Sec. 5]{step}).
The mixed derivative of $\Omega_m(\xb)$ is
$$\omega_m(\xb)=1-\frac{4}{m}\sum_{j}x_j +\frac{8}{m(m-1)}\sum_{i<j}x_i x_j$$
and
$$\int_{I^m}\omega_m^2(\xb)\,d\xb =\frac{2}{9m(m-1)},\quad \int_{I^m}\Omega_m(\xb)d\xb=\int_{I^m} \prod_{j=1}^m x_j \,\omega_m(\xb)\, d\xb=\frac{1}{9\cdot 2^m }.$$
Therefore, according to (\ref{absS}) and (\ref{absW})
\begin{gather}\label{eff}
e_S(\Omega_m)=e_W(\Omega_m)=\frac{m(m-1)}{18((4/3)^m -(m/3)-1)}.
\end{gather}
Again, the test statistics $S_{m,n}$ and $W_{m,n}$ are equally efficient in the Pitman sense. Their asymptotic efficiency
decreases as $m$ increases, and equals $ 0.9000$, $0.8060,$ $0.7181$ for  $m=3,\,4,\,5,$ respectively.


In both examples the asymptotic equivalence (in the sense of Pitman) of the tests based on $S_{m,n}$ and $W_{m,n}$ is explained by the fact
 that the corresponding cdfs in model (\ref{mod}) are radially symmetric, i.e., $F_{\t}(\xb)=\bar{F}_{\t}(1-\xb),$
 in which case $s_m(F_\t)$ and $w_m(F_\t)$ are known to be equal (see \cite[Sec. 3]{s-s}).
\bigskip

\n{\bf\large Acknowledgments}

\bigskip
The research of A. Nazarov was partly supported by RFBR grant 07-01-00159.
The research of N. Stepanova was supported by an NSERC grant.
We would like to thank Prof. Yu. V. Tarannikov for communicating us references
\cite{kl} and \cite{korsh}.

\end{document}